# ON OPTIMALITY OF STEPDOWN AND STEPUP MULTIPLE TEST PROCEDURES

By E. L. Lehmann, Joseph P. Romano and Juliet Popper Shaffer

*University of California, Berkeley, Stanford University and University of California, Berkeley*

Consider the multiple testing problem of testing $k$ null hypotheses, where the unknown family of distributions is assumed to satisfy a certain monotonicity assumption. Attention is restricted to procedures that control the familywise error rate in the strong sense and which satisfy a monotonicity condition. Under these assumptions, we prove certain maximin optimality results for some well-known stepdown and stepup procedures.

**1. Introduction.** For classical single-stage multiple comparison procedures, a number of optimality results are available. (See, e.g., [6], Chapter 11, and [11], Chapter 7 particularly Sections 7.9 and 7.10.) However, no such literature exists for the more recent stepdown and stepup methods. It is the purpose of the present paper to establish optimality properties for procedures of the latter kind.

Our setup and conditions are those of Lehmann [9], who discusses such an optimality result for the testing of two hypotheses. For the general problem of testing $k$ null hypotheses $H_1, \dots, H_k$, consider $k$ random variables $X_1, \dots, X_k$; typically, these are test statistics for the individual hypotheses $H_1, \dots, H_k$. We assume that $X = (X_1, \dots, X_k)$ has some $k$-dimensional joint cumulative distribution function $F_\theta(\cdot)$ indexed by $\theta = (\theta_1, \dots, \theta_k)$ in $\mathbf{R}^k$. The null hypothesis $H_i$ states $\theta_i \le 0$, which is being tested against the alternatives $\theta_i > 0$.

Stepdown procedures were initiated by Holm [7], while the stepup approach can be found in [2, 5, 8, 14, 16]. Background material on stepwise procedures is provided by Hochberg and Tamhane [6] and Westfall and Young [18]. Roughly speaking, stepdown procedures start by rejecting the most significant hypothesis (corresponding to the largest $X_i$) and then they sequentially consider the most significant of the remaining hypotheses. Alternatively, stepup procedures start with the least significant hypothesis (corresponding to the smallest $X_i$).

Our optimality results require crucial monotonicity assumptions and restrictions. We say that a region $\mathcal{M}$ of $x$ values is a monotone increasing region if $x = (x_1, \dots, x_k) \in \mathcal{M}$ and $x_i \le y_i$ implies that $y = (y_1, \dots, y_k)$ is also in $\mathcal{M}$. We assume of our model that increased values of $\theta$ lead to higher values of $X$, specifically, that if $\theta_i \le \gamma_i$, then

$$\int_{\mathcal{M}} dF_\theta(x_1, \dots, x_k) \le \int_{\mathcal{M}} dF_\gamma(x_1, \dots, x_k) \tag{1}$$

for every monotone increasing region $\mathcal{M}$. This assumption holds, in particular, if the distributions $F_\theta$ have densities $p_\theta$ with (increasing) monotone likelihood ratio; that is, if $x = (x_1, \dots, x_k)$, $y = (y_1, \dots, y_k)$, $\theta = (\theta_1, \dots, \theta_k)$ and $\theta' = (\theta'_1, \dots, \theta'_k)$, then

$$\frac{p_{\theta'}(x)}{p_\theta(x)} \le \frac{p_{\theta'}(y)}{p_\theta(y)}$$

whenever $x_i \le y_i$ for all $i$ and $\theta_j \le \theta'_j$ for all $j$. This notion of monotonicity was studied in [10]; other notions of stochastic ordering are discussed in [12].

In addition to condition (1), we will assume an analogous monotonicity property for the distribution of $(\delta_1 X_1, \dots, \delta_k X_k)$, for any $\delta_i \in \{-1, 1\}$. Specifically, for every monotone increasing region $\mathcal{M}$ and $\delta_i\theta_i \le \delta_i\gamma_i$,

$$\int_{\mathcal{M}} dF_\theta(\delta_1 x_1, \dots, \delta_k x_k) \le \int_{\mathcal{M}} dF_\gamma(\delta_1 x_1, \dots, \delta_k x_k). \tag{2}$$

For example, the condition for $(-X_1, \dots, -X_k)$ means that for any monotone decreasing region $\mathcal{M}'$ (the complement of a monotone increasing region), the inequality (1) is reversed; that is, the probability of the event $\{X \in \mathcal{M}'\}$ increases as $\theta$ decreases (in each component).

Under these assumptions, we shall restrict attention to decision rules satisfying the following monotonicity condition. A decision rule $D$ based on $X$ states for each possible value $x$ of $X$ the subset $I = I_x$ of $\{1, \dots, k\}$ of values $i$ for which the hypothesis $H_i$ is rejected. A decision rule $D$ is said to be monotone if

$$x_i \le y_i \quad \text{for } i \in I_x \quad \text{but } y_i < x_i \quad \text{for } i \notin I_x$$

implies that $I_x = I_y$. Thus, the subset of $x$ values that results in rejecting all hypotheses is a monotone increasing region. More generally, fix $I \subset \{1, \dots, k\}$ and, based on a monotone decision rule, let $\mathcal{M}_I$ denote the set of $x$ values such that $I_x = I$. If $\delta_i = 1$ for $i \in I$ and $\delta_i = -1$ otherwise, then

$$\{(\delta_1 x_1, \dots, \delta_k x_k) : (x_1, \dots, x_k) \in \mathcal{M}_I\}$$

is a monotone increasing set. By assumption (2), the probability of this set is increasing in $\delta_i\theta_i$.

Among all monotone decision rules that provide strong control of the familywise error rate (FWER), that is, of the probability of committing a Type 1 error by wrongly rejecting one or more true hypotheses, under any configuration of true and false null hypotheses we shall show how to maximize certain aspects of the power of the procedures (i.e., of the probability of correctly rejecting false hypotheses). However, we note that we are not restricting attention to any kind of stepwise procedure; rather, the resulting optimal procedures take the form of well-known stepwise procedures, which will be fully described later.

Here the restriction to monotone procedures is not just for convenience—the results are not true without this restriction. It is, in fact, possible to improve the rejection probability without violating the error control by adding small implausible pieces to the rejection regions, resulting in decision rules that are very counterintuitive. That this is possible is due to the fact that the bound for the error control is not attained but only approached in the limit as some parameter values tend to $\infty$ or $-\infty$. For a discussion of the pros and cons of such counterintuitive decision rules with references to the literature, see [13].

To conclude this introduction, we mention some situations in which the present approach does and some in which it does not apply. As a first example, consider a paired comparison experiment with pairs of observations $(Y_i, Z_i)$. Let $E(Y_i) = \mu_i$ and $E(Z_i) = \nu_i$, and consider testing the hypotheses $\theta_i = \nu_i - \mu_i = 0$ against the alternatives $\theta_i > 0$. If we put $X_i = Z_i - Y_i$ and base our inferences on the $X$'s, this reduces to the situation considered here. This example can be extended to the comparison of two treatments with $m_i$ and $n_i$ observations $(i = 1, \dots, k)$, respectively, on $k$ subjects. Another application is the comparison of $k$ treatments with a control. Here $\theta_i = \mu_i - \mu_0$, where the $\mu_i$ $(i = 1, \dots, k)$ and $\mu_0$ are the means for the $k$ treatments and the control, respectively.

On the other hand, the approach does not apply to the comparison of $k$ treatments, that is, the hypothesis $H : \mu_1 = \cdots = \mu_k$, where in the case of rejection one wishes to determine the pairs $i < j$ for which $\mu_i < \mu_j$. As in the preceding examples, the hypothesis can be reduced to $H : \theta_2 = \cdots = \theta_k = 0$ with, for example, $\theta_i = \mu_i - \mu_1$. However, with the resulting procedure, one can only determine the significant differences $\mu_j - \mu_i$ with $i = 1$ and not those with $1 < i < j$.

In Section 2 we treat the case $k = 2$ separately. In Section 3 we consider general $k$ for stepdown procedures, but make a further exchangeability assumption. The corresponding results for stepup procedures are then provided in Section 4, though a further assumption of monotonicity of critical values is invoked. Section 5 is a brief conclusion and all proofs are deferred to Section 6.

*Distributional assumptions.* We suppose $(X_1, \ldots, X_k)$ has a joint cumulative distribution function $F_\theta(\cdot)$, indexed by $\theta = (\theta_1, \ldots, \theta_k)$ in $\mathbf{R}^k$. The parameter space is a finite or infinite open rectangle with

$$\theta_i^L < \theta_i < \theta_i^U, \qquad i = 1, \ldots, k.$$

Similarly, the sample space is assumed to be a finite or infinite open rectangle with

$$x_i^L < X_i < x_i^U,$$

independent of $\theta$. We further assume the distribution of any subcollection $\{X_i, i \in I\}$ depends only on those $\theta_i$ with $i \in I$, and that $X_i$ tends in probability to $x_i^L$ as $\theta_i \to \theta_i^L$ and $X_i$ tends in probability to $x_i^U$ as $\theta_i \to \theta_i^U$.

To ease the notation, we assume here and in the remainder of the paper that $\theta_i$ varies in all of $\mathbf{R}$, so that $\theta_i^L = -\infty$ and $\theta_i^U = \infty$. We also simplify the notation by taking $x_i^L = -\infty$ and $x_i^U = \infty$. In addition, we assume that the joint distribution of $X$ has a density with respect to Lebesgue measure; this is used only so that the critical constants of the optimal procedures can be obtained for control at a given level $\alpha$ to be achieved exactly, but this hypothesis can certainly be weakened. In order for the critical constants to be uniquely defined, we further assume that the joint density is positive on its (assumed rectangular) region of support, but this can be weakened as well.

## 2. The case $k = 2$.

We are testing hypotheses $H_1$ and $H_2$ with $H_i$ corresponding to $\theta_i \le 0$. Let $\omega_{0,0}$ denote the part of the parameter space where both $H_1$ and $H_2$ are true; let $\omega_{0,1}$ correspond to the part where $H_1$ is true and $H_2$ is not true; similarly for $\omega_{1,0}$ and $\omega_{1,1}$.

A decision rule $D$ analogously divides up the sample space into regions $d_{0,0}$, $d_{0,1}$, $d_{1,0}$ and $d_{1,1}$. For example, $d_{0,1}$ corresponds to the region in the sample space where $H_1$ is declared true and $H_2$ is declared false. Also, let $d_i$ be the region where $H_i$ is rejected, so $d_1 = d_{1,0} \cup d_{1,1}$ and $d_2 = d_{0,1} \cup d_{1,1}$.

We will restrict attention to rules $D$ that are

$$\text{monotone} \tag{3}$$

and such that the

$$\text{FWER} \le \alpha. \tag{4}$$

For $\varepsilon = (\varepsilon_1, \varepsilon_2)$ with $\varepsilon_i > 0$, consider subsets of elements $(\theta_1, \theta_2)$ defined by

$$A_1(\varepsilon) = \{\theta_1 > \varepsilon_1\} \cup \{\theta_2 > \varepsilon_2\} \tag{5}$$

and

$$A_2(\varepsilon) = \{\theta_1 > \varepsilon_1\} \cap \{\theta_2 > \varepsilon_2\}. \tag{6}$$

A decision rule is deemed good if the quantities

$$\inf_{\theta \in A_1(\varepsilon)} P_\theta\{d_{0,0}^c\} \tag{7}$$

and

$$\inf_{\theta \in A_2(\varepsilon)} P_\theta\{d_{1,1}\} \tag{8}$$

are large. As we will see, it is *not* possible to find a rule $D$ satisfying (3) and (4) that maximizes (7) and (8) simultaneously.

In order to appreciate the criteria (7) and (8), first suppose $\theta \in A_1(\varepsilon)$. Then at least one $\theta_i$ is positive and so we would not want to conclude that both $\theta_i$ are $\leq 0$; rather, we wish to conclude $d_{0,0}^c$. Thus, maximizing (7) maximizes the minimum probability that we do not conclude $d_{0,0}$ as $\theta$ varies in $A_1(\varepsilon)$. Similarly, if $\theta \in A_2(\varepsilon)$, then both $\theta_i$ are positive, and so we wish to maximize the (minimum) chance that we make the decision $d_{1,1}$.

In addition, we also consider the following notion of optimality. Again suppose $\theta \in A_1(\varepsilon)$, so that at least one $\theta_i$ is positive. Then, as above, we do not want to make the decision $d_{0,0}$. However, we also do not wish to make the decision $d_{0,1}$ if, in fact, $H_1$ is false and $H_2$ is true; we would rather make the correct decision $d_{1,0}$. So, we also consider the probability of maximizing

$$\inf_{\theta \in A_1(\varepsilon)} P_\theta\{\text{reject at least one false hypothesis}\}. \tag{9}$$

In other words, the criterion (7) maximizes the minimum probability of rejecting exactly least one hypothesis (regardless of which are true and false), while criterion (9) maximizes the minimum probability of rejecting at least one *false* hypothesis. The latter criterion seems more compelling, though the former criterion might be justified in a situation where it is important to know that the joint null hypothesis (i.e., the global hypothesis that both hypotheses are true) is not true. In any case, we shall see that the same optimal procedure $D$ arises from both criteria.

THEOREM 2.1. *Consider the case $k = 2$ under the assumptions given in Section 1.*

(i) *A rule $D$ satisfying (3) and (4) maximizes (7) if*

$$d_{0,0}^c = \{X_1 > a_1 \text{ or } X_2 > a_2\} \tag{10}$$

*and $\{X_i > a_i\} \subset d_i$, where $a_1$ and $a_2$ are determined so that*

$$P_{0,0}\{X_1 > a_1 \text{ or } X_2 > a_2\} = \alpha \tag{11}$$

*and*

$$P_{\varepsilon_1}\{X_1 > a_1\} = P_{\varepsilon_2}\{X_2 > a_2\}. \tag{12}$$

*Its minimum (rejection) probability over $A_1(\varepsilon)$ is given by $P_{\varepsilon_1}\{X_1 > a_1\}$.*

(ii) *The (stepdown) rule $D$ satisfying* (3), (4) *and* (10) *that maximizes* (8) *is given by*

$$d_{0,1} = \{X_1 < b_1, X_2 \geq a_2\}, \tag{13}$$

$$d_{1,0} = \{X_1 \geq a_1, X_2 < b_2\}, \tag{14}$$

$$d_{1,1} = \{X_1 \geq b_1, X_2 \geq b_2\} \cap d_{0,0}^c, \tag{15}$$

*where $b_i$ satisfies*

$$P_0\{X_i \geq b_i\} = \alpha \tag{16}$$

*(and so $b_i < a_i$). The minimum probability of $d_{1,1}$ over $A_2(\varepsilon)$ is*

$$P_{\varepsilon_1,\varepsilon_2}\{X_1 > a_1, X_2 > b_2 \cup X_1 > b_1, X_2 > a_2\}.$$

(iii) *The result* (i) *holds for $D$ if criterion* (7) *is replaced by* (9), *and* (12) *is also the maximum value of criterion* (9).

Note that once $d_{0,0}$ and $d_{1,1}$ are determined, so are $d_{0,1}$ and $d_{1,0}$ (by monotonicity).

The procedure $D$ of Theorem 2.1 is an example of a stepdown procedure. It starts by rejecting the most significant hypothesis (corresponding to the largest $X_i$) and it then sequentially considers the most significant of the remaining hypotheses. Alternatively, stepup procedures start with the least significant hypothesis (corresponding to the smallest $X_i$), and an optimality result is now given for such a procedure.

REMARK 2.1. The proof shows that the optimal procedure $D$ in (i) and (ii) is the unique rule satisfying (3) and (4) which maximizes (7), in the sense that if $E$ is any other such rule, then $e_{0,0} \triangle d_{0,0}$ has Lebesgue measure 0, where $A \triangle B$ denotes the symmetric difference between sets $A$ and $B$. Similarly, a rule $E$ satisfying (3), (4) and (10) maximizing (8) must satisfy $e_{1,1} \triangle d_{1,1}$ has Lebesgue measure 0.

Also, notice that the optimal procedure $D$ does not depend on $\varepsilon$. It follows that $D$ is admissible in the following sense. Suppose there exists another monotone rule $E$ that controls the FWER, and such that

$$P_\theta\{d_{0,0}^c\} \leq P_\theta\{e_{0,0}^c\} \qquad \text{for all } \theta \in \omega_{0,0}^c, \tag{17}$$

with strict inequality for some $\theta \in \omega_{0,0}^c$. Taking the infimum of both sides over $\theta \in A_1(0)$, it follows that $E$ must also be optimal in the sense of Theorem 2.1(i). But, by uniqueness, $e_{0,0} \triangle d_{0,0}$ has Lebesgue measure 0, which implies the $\leq$ in (17) is an equality. A similar admissibility result for the region $d_{1,1}$ can be stated as well.

Analogous uniqueness and admissibility results hold for all the optimal procedures presented later on. For a discussion of admissibility in multiple testing problems, see [3].

THEOREM 2.2. *Consider the case* $k = 2$ *under the assumptions given in Section* 1.

(i) *A rule* $D$ *satisfying* (3) *and* (4) *maximizes* (8) *if* $d_{1,1}$ *is given by*

$$d_{1,1} = \{X_1 > b_1, X_2 > b_2\}, \tag{18}$$

*and* $d_i \subset \{X_i > b_i\}$, *where* $b_i$ *satisfies* (16) *(so it is the same constant as in Theorem* 2.1*). Its minimum probability over* $A_2(\varepsilon)$ *is given by* $P_{\varepsilon_1,\varepsilon_2}\{X_1 > b_1, X_2 > b_2\}$.

(ii) *The (stepup) rule* $D$ *satisfying* (3), (4) *and* (18) *that maximizes* (7) *is given by*

$$d_{0,1} = \{X_1 < b_1, X_2 \geq \tilde{a}_1\}, \tag{19}$$

$$d_{1,0} = \{X_1 \geq \tilde{a}_1, X_2 < b_2\}, \tag{20}$$

$$d_{0,0} = \{X_1 \leq \tilde{a}_1, X_2 \leq \tilde{a}_2\} \cap d_{1,1}^c, \tag{21}$$

*where* $\tilde{a}_i$ *is determined so that*

$$P_{0,0}\{d_{0,0}^c\} = \alpha \tag{22}$$

*and*

$$P_{\varepsilon_1}\{X_1 \geq \tilde{a}_1\} = P_{\varepsilon_2}\{X_2 \geq \tilde{a}_2\}. \tag{23}$$

*The value of* (23) *is the minimum probability of* $D$ *over* $A_1(\varepsilon)$.

(iii) *The result* (ii) *holds for* $D$ *if criterion* (7) *is replaced by* (9), *and* (23) *is also the maximum value of criterion* (9).

REMARK 2.2. Note that $b_i < a_i < \tilde{a}_i$. Also, the best minimum probability over $A_1(\varepsilon)$ in the case of Theorem 2.1 exceeds the best in the case of Theorem 2.2, but it reverses for Theorem 2.2.

REMARK 2.3. A remark analogous to Remark 2.1 applies to the optimal procedure in Theorem 2.2.

REMARK 2.4. It is now clear that, subject to (3) and (4), we cannot find a rule to maximize both (7) and (8). By Theorem 2.1(i) and Theorem 2.2(ii), such a rule $D$ would have to satisfy

$$d_{0,0} = \{X_1 \leq a_1 \text{ and } X_2 \leq a_2\}$$

and

$$d_{1,1} = \{X_1 > b_1, X_2 > b_2\}$$

simultaneously, which is impossible because these two sets have a nontrivial intersection as $b_i < a_i$.

REMARK 2.5. The results of this paper do not hold without the monotonicity assumption. For example, consider part (i) of Theorem 2.2. Suppose further that $X_1$ and $X_2$ are independent with $X_i$ normally distributed with mean $\theta_i$ and variance 1. Then $b_i = b = z_{1-\alpha}$, the $1-\alpha$ quantile of the standard normal distribution. The probability of $d_{1,1}$ under $(\theta_1, \theta_2)$ with both $\theta_i > 0$ is always less than $\alpha$ and approaches $\alpha$ as either $\theta_i \to \infty$. Therefore, by adding to $d_{1,1}$ a small enough region in the southwest quadrant, one can increase the rejection probability without violating the level constraint; see Section 4 of [13]. Such a procedure is not monotone. Similarly, regarding the problem addressed in (i) of Theorem 2.1, [9], Section 3, shows that the maximin test is not monotone.

**3. General $k$ stepdown.** Consider testing $k$ null hypotheses $H_1, \dots, H_k$ with $H_i$ corresponding to $\theta_i \le 0$. In this section and the next, we add a symmetry condition for the joint distribution of $(X_1, \dots, X_k)$. Specifically, we assume that the joint distribution of $(X_1, \dots, X_k)$ under $\theta_i = \theta$ (some value independent of $i$) is exchangeable. This is not a crucial assumption (and actually only needs to hold at $\theta = 0$ or $\theta = \varepsilon$, where $\varepsilon$ is given in the statement of the theorems), but it reduces the number of critical values from order $2^k$ to $k$. The results should generalize, but at the expense of more complicated notation.

Let

$$X_{r_1} \ge X_{r_2} \ge \cdots \ge X_{r_k}$$

denote the ordered $X$-values, and let $H_{r_1}, \dots, H_{r_k}$ denote their corresponding null hypotheses.

For any (monotone) decision rule $E$, let $E_{k,j}$ denote the event that $E$ rejects at least $j$ of the null hypotheses. For $\varepsilon > 0$, let

$$A_j(\varepsilon) = \{(\theta_1, \dots, \theta_k) : \text{at least } j \ \theta_i \text{ satisfy } \theta_i > \varepsilon\}.$$

Consider the monotone stepdown decision rule $D$ that rejects $H_{r_1}, \dots, H_{r_j}$ and accepts the remaining null hypotheses if $X_{r_i} \ge c_{k,i}$ for $1 \le i \le j$, but $X_{r_j+1} < c_{k,j+1}$, where the $c_{k,j} = c_{k,j}(\alpha)$ are determined by

$$P_{\underbrace{0,\dots,0}_{k-j+1 \text{ times}}} \{X_i > c_{k,j} \text{ for some } i, 1 \le i \le k-j+1\} = \alpha. \tag{24}$$

Then

$$D_{k,j} = \{X_{r_i} \ge c_{k,i}, 1 \le i \le j\}.$$

Note that (24) implies the important relationship

$$c_{k,j} = c_{k-1,j-1} \tag{25}$$

if $k \geq j \geq 2$. Also note that, for fixed $k$, $c_{k,j}$ is nonincreasing in $j$.

Since the constants $c_{k,j}$ depend only on $k-j$, we may more succinctly define

$$f_{k-i+1} \equiv c_{k,i}, \tag{26}$$

where the $f_j$ are determined by

$$P_{\underbrace{0,\ldots,0}_{j \text{ times}}}\{\max(X_1,\ldots,X_j) > f_j\} = \alpha. \tag{27}$$

The procedure $D$ then rejects $H_{r_1},\ldots,H_{r_j}$ if and only if $X_{r_i} \geq f_{k-i+1}$ for $1 \leq i \leq j$.

LEMMA 3.1. *Suppose the assumptions of Section 1 and the symmetry condition described at the beginning of this section hold.*

(i) *The above decision rule $D$ controls the FWER at level $\alpha$.*
(ii) *Define*

$$\beta_{k,j}(\alpha,\varepsilon) = \inf_{\theta \in A_j(\varepsilon)} P_\theta\{D_{k,j}\};$$

*that is, $\beta_{k,j}(\alpha,\varepsilon)$ is the minimum probability of $D_{k,j}$ over $A_j(\varepsilon)$. Then*

$$\beta_{k,j}(\alpha,\varepsilon) = P_{\underbrace{\varepsilon,\ldots,\varepsilon}_{j \text{ times}}}\{S_{k,j}\}, \tag{28}$$

*where*

$$\begin{aligned} S_{k,j} = \{&X_{\pi_j(1)} > f_k,\ldots,X_{\pi_j(j)} > f_{k-j+1} \\ &\text{for some permutation } \pi_j \text{ of } \{1,\ldots,j\}\}. \end{aligned} \tag{29}$$

*So (28) is the minimum probability over $A_j(\varepsilon)$ not only of rejecting at least $j$ hypotheses, but also of rejecting at least $j$ false hypotheses.*

THEOREM 3.1. *Suppose the assumptions of Section 1 and the symmetry condition described at the beginning of this section hold.*

(i) *Among monotone decision rules $E$ that control the FWER, $D$ maximizes*

$$\inf_{\theta \in A_1(\varepsilon)} P_\theta\{E_{k,1}\}. \tag{30}$$

*Also, $D$ maximizes*

$$\inf_{\theta \in A_2(\varepsilon)} P_\theta\{E_{k,2}\}$$

*among such rules $E$ that also satisfy $E_{k,2} \subset D_{k,1}$. In general, for $j = 2, \dots, k$, $D$ maximizes*

$$\inf_{\theta \in A_j(\varepsilon)} P_\theta\{E_{k,j}\} \tag{31}$$

*among monotone rules $E$ that control the FWER and satisfy*

$$E_{k,j} \subset D_{k,j-1}. \tag{32}$$

*Therefore, for any other rule $E$, we must have*

$$\inf_{\theta \in A_j(\varepsilon)} P_\theta(E_{k,j}) < \beta_{k,j}(\alpha, \varepsilon)$$

*for at least one $j$.*

(ii) *$D$ also is optimal in the sense that it maximizes*

$$\inf_{\theta \in A_j(\varepsilon)} P_\theta\{\textit{reject at least } j \textit{ false hypotheses}\}$$

*subject to* (32).

Remark 3.1. The procedure $D$ is essentially unique (up to sets of Lebesgue measure 0), as described in Remark 2.1, and an admissibility result analogous to that described in Remark 2.1 holds as well.

Remark 3.2. For fixed $k$, the optimal constants with $c_{k,j} = f_{k-j+1}$ are given by the values

$$c_{k,1}, c_{k,2}, \dots, c_{k,k}. \tag{33}$$

But, since $c_{k,2} = c_{k-1,1}$, $c_{k,3} = c_{k-2,1}$, and so on, the sequence (33) is equivalent to

$$c_{k,1}, c_{k-1,1}, c_{k-2,1}, \dots, c_{1,1}.$$

This is just a sequentially rejective procedure of the kind proposed by Holm [7]: after the first step using the critical value $c_{k,1}$, reduce the number of hypotheses from $k$ to $k-1$ and repeat the first step but now with $c_{k-1,1}$, and so on. In the case where the $X_i$ have a uniform $(0,1)$ marginal distribution under the null hypothesis so that we translate everything into $p$-values and reject for *small* values, Holm [7] used $c_{k,1} = \alpha/k$ since he assumed only the marginal distributions to be known (and strong error control follows by Bonferroni). Our $c_{k,1}$ would then be determined by

$$P_{\underbrace{0,\dots,0}_{k \text{ times}}}\{X_i \le c_{k,1} \text{ for one or more values of } i : 1 \le i \le k\} = \alpha$$

or, equivalently,

$$P_{\underbrace{0,\dots,0}_{k \text{ times}}}\{\min(X_1, \dots, X_k) \le c_{k,1}\} = \alpha.$$

If we further assume independence of the $p$-values, then the critical constants $c_{k,j}$ satisfy

$$1-(1-c_{k,j})^{k-j+1}=\alpha.$$

Thus, the Holm principle remains in effect, except that instead of using $c_{k,1}=\alpha/k$, the independence assumption implies the exact critical values $c_{k,1}=1-(1-\alpha)^{1/k}$.

**4. General *k* stepup.** Assume the conditions imposed in the previous section. We are testing null hypotheses $H_1,\ldots,H_k$ with $H_i$ corresponding to $\theta_i\le 0$. Let

$$X_{(1)}\le X_{(2)}\le\cdots\le X_{(k)}$$

denote the ordered $X$-values; in the notation of the previous section, $X_{(j)}=X_{r_{k-j+1}}$.

Consider the following monotone stepup decision rule $D$ for appropriately chosen constants $d_1,\ldots,d_k$ (to be specified shortly, but assumed nondecreasing). If $X_{(1)}>d_1$, then reject all null hypotheses. Otherwise, if $X_{(1)}\le d_1$ but $X_{(2)}>d_2$, reject the $k-1$ hypotheses corresponding to the $k-1$ largest $X$'s. In general, for the smallest $j$ such that $X_{(j)}>d_j$, reject the $k-j+1$ hypotheses corresponding to the $k-j+1$ largest $X$'s and accept the remaining. (Note that the constants $d_j$ should perhaps be written as $d_{k,j}$ to show the dependence on $k$; however, we will see that $d_{k,j}$ will be chosen to be independent of $k$ and so we just abbreviate to $d_j$.)

The above rule rejects at least $j$ null hypotheses for the set $D_{k,j}$ defined by

$$D_{k,j}=\{X_{(1)}>d_1\}\cup\cdots\cup\{X_{(k-j+1)}>d_{k-j+1}\}.$$

Equivalently, at least $k-j+1$ hypotheses are accepted if $D^c_{k,j}$ occurs, where

$$D^c_{k,j}=\{X_{(1)}\le d_1\}\cap\cdots\cap\{X_{(k-j+1)}\le d_{k-j+1}\}.$$

Evidently,

$$D_{k,j+1}\subset D_{k,j}.$$

The constants $d_j$ are determined so that

$$P_{\underbrace{0,\ldots,0}_{j\text{ times}}}\{L_j\}=1-\alpha, \tag{34}$$

where

$$L_j=\{X_{\pi(1)}\le d_1,\ldots,X_{\pi(j)}<d_j\text{ for some permutation of }\{1,\ldots,j\}\}. \tag{35}$$

Note that the constant $d_j$ does not depend on $k$ as reflected in the notation. Also, $d_1=c_{k,k}=c_{1,1}=f_1$, where $c_{1,1}$ and $f_1$ are the constants (24) and (26) of the previous section. However, as pointed out by an anonymous referee, in the case $k > 2$, it need not be the case that $d_{k,j}$ is nondecreasing in $j$. A counterexample is provided in [4]; some further references on the monotonicity of critical values are [1] and [15]. In order to prove our results, we need to assume the monotonicity holds.

LEMMA 4.1. *Assume the conditions of Lemma 3.1. Also assume that the constants $d_j$ used in the procedure $D$ are nondecreasing in $j$.*

(i) *The above decision rule $D$ controls the FWER at level $\alpha$.*
(ii) *Define*

$$\tilde{\beta}_{k,j} = \inf_{\theta \in A_j(\varepsilon)} P_\theta\{D_{k,j}\};$$

*that is, $\tilde{\beta}_{k,j}(\alpha, \varepsilon)$ is the minimum probability of $D_{k,j}$ over $A_j(\varepsilon)$. Then*

$$\tilde{\beta}_{k,j}(\alpha, \varepsilon) = P_{\underbrace{\varepsilon,\ldots,\varepsilon}_{j\ \text{times}}}\{\min(X_1, \ldots, X_j) > d_{k-j+1}\}. \tag{36}$$

*The minimum probability over $A_j(\varepsilon)$ of rejecting at least $j$ false hypotheses is also given by (36).*

THEOREM 4.1. *Assume the conditions of Theorem 3.1. Also assume that the constants $d_j$ used in the procedure $D$ are nondecreasing in $j$.*

(i) *Among monotone decision rules $E$ that control the FWER at level $\alpha$, $D$ maximizes*

$$\inf_{\theta \in A_k(\varepsilon)} P_\theta\{E_{k,k}\}. \tag{37}$$

*Also, $D$ maximizes*

$$\inf_{\theta \in A_{k-1}(\varepsilon)} P_\theta\{E_{k,k-1}\}$$

*among rules that satisfy $D_{k,k} \subset E_{k,k-1}$. In general, for $j = k-1, \ldots, 1$, $D$ maximizes*

$$\inf_{\theta \in A_j(\varepsilon)} P_\theta\{E_{k,j}\} \tag{38}$$

*among monotone rules $E$ that control the FWER and satisfy*

$$D_{k,j+1} \subset E_{k,j}. \tag{39}$$

*Therefore, for any other rule $E$, we must have*

$$\inf_{\theta \in A_j(\varepsilon)} P_\theta(E_{k,j}) < \tilde{\beta}_{k,j}(\alpha, \varepsilon)$$

*for at least one $j$.*

(ii) *$D$ also is optimal in the sense that it maximizes*

$$\inf_{\theta \in A_j(\varepsilon)} P_\theta\{\textit{reject at least } j \textit{ false hypotheses}\}$$

*subject to* (39).

Remark 4.1. Again, the procedure $D$ is unique up to sets of Lebesgue measure 0, and it is admissible; see Remark 2.1.

Remark 4.2. Letting

$$X_{j:1} \leq \cdots \leq X_{j:j}$$

denote the ordered values of just the first $j$ $X$'s, the constants $d_j$ are determined by

$$P_{\underbrace{0,\ldots,0}_{j \text{ times}}}\{X_{j:1} \leq d_1, \ldots, X_{j:j} \leq d_j\} = 1 - \alpha.$$

If we compare this with (27), we see that $f_j < d_j$, except when $j = 1$, in which case $f_1 = d_1$.

**5. Conclusions.** Stepdown and stepup methods were proposed as intuitively appealing by Holm, Hochberg, Dunnett and Tamhane, and others. The present paper, treating the case of one-sided alternatives only, used optimality criteria that seemed reasonable and were not selected to justify predetermined solutions. It is gratifying that the results confirm the intuition of the originators of these methods. Even though our assumptions are strong, some stepwise methods can now be viewed as asymptotically optimal, such as the stepup method of Dunnett and Tamhane [2]. Outside the strong assumptions imposed in this paper, Westfall and Young [18] give general resampling methods to approximate the critical values of stepdown procedures, while Troendle [17] addresses the corresponding problem for stepup procedures.

**6. Proofs and auxiliary results.**

Proof of Theorem 2.1. First, observe that for the procedure $D$ given in (i), for $\theta \in \omega_{0,0}$,

$$P_\theta\{d_{0,0}^c\} \leq P_{0,0}\{d_{0,0}^c\} = P_{0,0}\{X_1 > a_1 \text{ or } X_2 > a_2\} = \alpha$$

by choice of $a_i$. For this $D$, by monotonicity, the inf over $\theta \in A_1(\varepsilon)$ in (7) occurs at $(\theta_1, \theta_2) = (\varepsilon_1, -\infty)$ or $(-\infty, \varepsilon_2)$; this is a shorthand notation so that

$$P_{\varepsilon_1,-\infty}\{d_{0,0}^c\} = \lim_{\theta_2 \to -\infty}\{d_{0,0}^c\}.$$

But then

$$P_{\varepsilon_1,-\infty}\{d_{0,0}^c\} = P_{\varepsilon_1}\{X_1 > a_1\}$$

and

$$P_{-\infty,\varepsilon_2}\{d_{0,0}^c\} = P_{\varepsilon_2}\{X_2 > a_2\}.$$

So, the value of criterion (7) for the procedure $D$ is indeed given by (12). Similarly, the value of criterion (9) for $D$ is also (12). Indeed, as $\theta_1 \to -\infty$, the chance that $H_1$ is incorrectly rejected tends to 0.

To prove (i), suppose $E$ is another decision rule satisfying (3) and (4). Assume there exists $(x_1, x_2) \notin d_{0,0}$, but $(x_1, x_2) \in e_{0,0}$. Then there exists at least one component with $x_i > a_i$, say $x_1 > a_1$. Hence,

$$P_{\varepsilon_1,-\infty}\{e_{0,0}) \ge P_{\varepsilon_1,-\infty}\{X_1 \le x_1, X_2 \le x_2\} = P_{\varepsilon_1}\{X_1 \le x_1\} > P_{\varepsilon_1}\{X_1 \le a_1\}.$$

Therefore,

$$P_{\varepsilon_1,-\infty}\{e_{0,0}^c\} < P_{\varepsilon_1,-\infty}\{X_1 \le a_1\} = P_{\varepsilon_1}\{X_1 \le a_1\},$$

so that $E$ has a smaller value of criterion (7) than does a claimed optimal $D$. So it must be the case that $e_{0,0} \subset d_{0,0}$. But, if $e_{0,0}$ is strictly contained in $d_{0,0}$ such that the set difference $e_{0,0} \triangle d_{0,0}$ has positive Lebesgue measure, then its region for rejecting $\omega_{0,0}$, namely, $e_{0,0}^c$, is bigger than $d_{0,0}^c$, implying

$$P_{0,0}\{e_{0,0}^c\} > P_{0,0}\{d_{0,0}^c\} = \alpha.$$

The conclusion is that an optimal region $D$ must have the stated region (10) $d_{0,0}^c$.

To prove (ii), let us first check that the claimed solution controls the FWER. For $\theta \in \omega_{0,0}$,

$$P_\theta\{d_{0,0}^c\} \le \alpha$$

as previously argued. For $\theta \in \omega_{0,1}$,

$$P_\theta\{\text{Type 1 error}\} = P_\theta\{d_{1,1} \cup d_{1,0}\} \le P_\theta\{X_1 \ge b_1\} \le P_0\{X_1 \ge b_1\} = \alpha$$

similarly for $\omega_{1,0}$.

The goal now is to find $D$ satisfying (3), (4) and $d_{0,0}^c$ given by (10) to maximize (8). Consider another rule $E$ satisfying (3), (4) and $e_{0,0} = d_{0,0}$. Suppose there exists $(x_1, x_2) \in e_{1,1}$ such that $x_i < b_i$ for some $i$, say $i = 1$. Then

$$P_{0,\infty}\{e_{1,1}\} \ge P_{0,\infty}\{X_1 \ge x_1, X_2 \ge x_2\} = P_0\{X_1 \ge x_1\} > P_0(X_1 \ge b_1) = \alpha,$$

which would contradict strong control. So $e_{1,1} \subset d_{1,1}$. But you cannot take away points from $d_{1,1}$ without lowering the minimum power at $(\theta_1, \theta_2) = (\varepsilon, \varepsilon)$.

To prove (iii), simply observe, for any $\theta$,

$$P_\theta\{\text{rejecting at least one false } H_i\} \le P_\theta\{\text{rejecting at least one } H_i\},$$

and so

$$\inf_{\theta\in A_1(\varepsilon)} P_\theta\{\text{rejecting at least one false } H_i\} \le \inf_{\theta\in A_1(\varepsilon)} P_\theta\{\text{rejecting at least one } H_i\}.$$

But the right-hand side is $P_{\varepsilon_1}\{X_1 > a_1\}$, and so it suffices to show that $D$ satisfies

$$\inf_{\theta\in A_1(\varepsilon)} P_\theta\{D \text{ rejects at least one false } H_i\} = P_{\varepsilon_1}\{X_1 > a_1\}.$$

But the earlier argument for (12) showed this to be the case. $\square$

Proof of Theorem 2.2. To prove (i), suppose $E$ is another rule satisfying (3) and (4) which rejects both hypotheses if $(X_1, X_2) \in e_{1,1}$. Suppose there exists $(x_1, x_2) \in e_{1,1}$ such that $x_i < b_i$ for some $i$, say $i = 1$. Then

$$P_{0,\infty}\{e_{1,1}\} \ge P_{0,\infty}\{X_1 \ge x_1, X_2 \ge x_2\} = P_0\{X_1 \ge x_1\} > P_0\{X_1 \ge b_1\} = \alpha,$$

which would contradict $E$ control of the FWER. So, $e_{1,1} \subset d_{1,1}$. But you cannot take away any point from $d_{1,1}$ without lowering the minimum power at $(\varepsilon_1, \varepsilon_2)$.

To prove (ii), note that, for the claimed solution the value of (7) is given by

$$\inf_{\theta\,:\,\theta\in A_1(\varepsilon)} P_\theta(d^c_{0,0}\} = P_{\varepsilon_1,-\infty}\{d^c_{0,0}\} = P_{\varepsilon_1}\{X_1 > \tilde{a}_1\}.$$

We now seek to determine $d_{0,0}$ [like Theorem 2.1(i) with the added constraint that $d_{0,0} \subset d^c_{1,1}$]. To prove optimality of the claimed solution, suppose $E$ is another rule satisfying (3), (4) and $e_{1,1} = d_{1,1}$, with $d_{1,1}$ given by (18). Suppose $(x_1, x_2) \notin d_{0,0}$, but $(x_1, x_2) \in e_{0,0}$, so that $x_i > \tilde{a}_i$ for some $i$, say $i = 1$. Then

$$\begin{aligned} P_{\varepsilon_1,-\infty}\{e_{0,0}\} &\ge P_{\varepsilon_1,-\infty}\{X_1 \le x_1, X_2 \le x_2\} \\ &= P_{\varepsilon_1}\{X_1 \le x_1\} > P_{\varepsilon_1}\{X_1 > \tilde{a}_1\}. \end{aligned}$$

Therefore,

$$P_{\varepsilon,-\infty}\{e^c_{0,0}\} < P_\varepsilon\{X_1 > \tilde{a}_1\},$$

so that $E$ cannot be optimal. So it must be the case that $e_{0,0} \subset d_{0,0}$. But if $e_{0,0}$ is strictly contained in $d_{0,0}$, its region for rejecting $\omega_{0,0}$, namely, $e^c_{0,0}$, is bigger than $d^c_{0,0}$, in which case

$$P_{0,0}\{e^c_{0,0}\} > P_{0,0}\{d^c_{0,0}\} = \alpha,$$

a contradiction of strong control.

Finally, we check that $D$ itself exhibits control of the FWER. For $\theta \in \omega_{0,0}$, the probability of a Type 1 error is $\leq \alpha$ because of (22). For $\theta = (\theta_1, \theta_2) \in \omega_{0,1}$,

$$\begin{aligned} P_\theta\{\text{Type 1 error}\} &\leq P_{0,\infty}\{X_1 > b_1, X_2 > b_2 \cup X_1 \geq \tilde{a}_1, X_2 < b_2\} \\ &= P_0\{X_1 > b_1\} = \alpha, \end{aligned}$$

as required.

The proof of (iii) is completely analogous to the proof of Theorem 2.1(iii). □

Proof of Lemma 3.1. To prove (i), suppose $H_1, \ldots, H_p$ are true and $H_{p+1}, \ldots, H_k$ are false. A Type 1 error occurs if any of $H_1, \ldots, H_p$ are rejected. For the rule $D$, the set where a rejection of any of $H_1, \ldots, H_p$ occurs is a monotone increasing set, and so the probability of this event is largest under this configuration of true and false hypotheses when

$$(\theta_1, \ldots, \theta_k) = (\underbrace{0, \ldots, 0}_{p \text{ times}}, \underbrace{\infty, \ldots, \infty}_{k-p \text{ times}}),$$

and this probability is equal to

$$P_{\underbrace{0,\ldots,0}_{p \text{ times}}}\{X_i > f_p \text{ for some } i = 1, \ldots, p\} = \alpha$$

by (27) with $j = p$.

To prove (ii), note that the minimum power occurs when $\theta$ is one of the $\binom{k}{j}$ points with $j$ values of $\varepsilon$ and $k - j$ values of $-\infty$, such as

$$w_{k,j} = w_{k,j}(\varepsilon) = (\underbrace{\varepsilon, \ldots, \varepsilon}_{j \text{ times}}, \underbrace{-\infty, \ldots, -\infty}_{k-j \text{ times}}). \tag{40}$$

Then, $P_{w_{k,j}}(D_{k,j})$ reduces to $\beta_{k,j}(\alpha, \varepsilon)$ as claimed. Also, for such a configuration $w_{k,j}$, only the $j$ hypotheses $H_1, \ldots, H_j$ can be rejected, and so the minimum probability of rejecting at least $j$ hypotheses is the same as the minimum probability of rejecting exactly $j$ hypotheses (and it is also equal to the probability of rejecting exactly $j$ false hypotheses). □

Before the proof of Theorem 3.1, we need two lemmas. We will make use of the following notation. If $R$ is any region in $\mathbf{R}^k$, let

$$R^z = \{(x_1, \ldots, x_{k-1}) : (x_1, \ldots, x_{k-1}, z) \in R\}.$$

Lemma 6.1. *Let $R$ be any monotone rejection region in $\mathbf{R}^k$ [so $x = (x_1, \ldots, x_k) \in R$ implies $y \in R$ if $y_i \geq x_i$ for all $i$].*

(i) *If* $z_1 < z_2$, *then* $R^{z_1} \subset R^{z_2}$.
(ii) $R^z$, $\bigcup_z R^z$ *and* $\bigcap_z R^z$ *are all monotone rejection regions in* $\mathbf{R}^{k-1}$.

Proof. If $z_1 < z_2$ and $(x_1, \ldots, x_{k-1}) \in R^{z_1}$, then $(x_1, \ldots, x_{k-1}, z_1) \in R$. By monotonicity, $(x_1, \ldots, x_{k-1}, z_2) \in R$, and so $(x_1, \ldots, x_{k-1}) \in R^{z_2}$. The proof of (ii) is just as easy. □

Lemma 6.2. *Assume the distributional assumptions given at the end of Section 1. Let R be any monotone rejection region in* $\mathbf{R}^k$. *Then for any values of the parameters* $\theta_1, \ldots, \theta_{k-1}$,

$$P_{\theta_1,\ldots,\theta_{k-1},\infty}(R) = P_{\theta_1,\ldots,\theta_{k-1}}\left\{\bigcup_z R^z\right\} \tag{41}$$

*and*

$$P_{\theta_1,\ldots,\theta_{k-1},-\infty}\{R\} = P_{\theta_1,\ldots,\theta_{k-1}}\left\{\bigcap_z R^z\right\}. \tag{42}$$

Proof. To prove (41),

$$\begin{aligned}P_{\theta_1,\ldots,\theta_{k-1},\infty}\{R\} &= \lim_{\theta_k\to\infty} P_{\theta_1,\ldots,\theta_k}\{(X_1, \ldots, X_{k-1}) \in R^{X_k}\}\\ &= \lim_{\theta_k\to\infty} P_{\theta_1,\ldots,\theta_k}\{(X_1, \ldots, X_{k-1}) \in R^{X_k}, X_k \ge z\}\\ &\le P_{\theta_1,\ldots,\theta_{k-1}}\left\{(X_1, \ldots, X_{k-1}) \in \bigcup R^z\right\}.\end{aligned}$$

Also, for every $z$,

$$\begin{aligned}P_{\theta_1,\ldots,\theta_{k-1},\infty}\{R\} &= \lim_{\theta_k\to\infty} P_{\theta_1,\ldots,\theta_k}\{(X_1, \ldots, X_{k-1}) \in R^{X_k}, X_k \ge z\}\\ &\ge P_{\theta_1,\ldots,\theta_{k-1}}\{(X_1, \ldots, X_{k-1}) \in R^z\}\end{aligned}$$

and so

$$P_{\theta_1,\ldots,\theta_{k-1},\infty}\{R\} \ge P_{\theta_1,\ldots,\theta_{k-1}}\left\{(X_1, \ldots, X_{k-1}) \in \bigcup R^z\right\},$$

and (41) follows.

To prove (42),

$$\begin{aligned}&P_{\theta_1,\ldots,\theta_{k-1},-\infty}(R)\\ &\quad= \lim_{\theta_k\to-\infty} P_{\theta_1,\ldots,\theta_k}\{(X_1, \ldots, X_{k-1}) \in R^{X_k}\}\\ &\quad= \lim_{\theta_k\to-\infty} P_{\theta_1,\ldots,\theta_k}\{(X_1, \ldots, X_{k-1}) \in R^{X_k} X_k, \le z\}\\ &\quad\le P_{\theta_1,\ldots,\theta_{k-1}}\{(X_1, \ldots, X_{k-1}) \in R^z\}\end{aligned}$$

for every $z$. Let $z \to -\infty$, so that $R^z$ decreases to $\bigcap R^z$. Then we can conclude

$$P_{\theta_1,\ldots,\theta_{k-1},-\infty}(R) \le P_{\theta_1,\ldots,\theta_{k-1}}\Big\{(X_1,\ldots,X_{k-1}) \in \bigcap R^z\Big\}.$$

Also,

$$\begin{aligned}P_{\theta_1,\ldots,\theta_{k-1},-\infty}\{R\} &= \lim_{\theta_k\to\infty} P_{\theta_1,\ldots,\theta_k}\{(X_1,\ldots,X_{k-1}) \in R^{X_k}\}\\ &\ge P_{\theta_1,\ldots,\theta_{k-1}}\Big\{(X_1,\ldots,X_{k-1}) \in \bigcap R^z\Big\},\end{aligned}$$

and the result follows. $\square$

Next, given a monotone rejection region $R$, define

$$U^1(R) = \bigcup_z R^z,$$

$$U^2(R) = U^1(U^1(R))$$

and

$$U^j(R) = U^1(U^{j-1}(R)).$$

Similarly, let

$$I^1(R) = \bigcap_z R^z$$

and

$$I^j(R) = I^1(I^{j-1}(R)).$$

By applying Lemma 6.2 repeatedly, we also obtain

$$P_{\theta_1,\ldots,\theta_{k-j},\underbrace{\infty,\ldots,\infty}_{j\ \text{times}}}\{R\} = P_{\theta_1,\ldots,\theta_{k-j}}\{U^j(R)\} \tag{43}$$

and

$$P_{\theta_1,\ldots,\theta_{k-j},\underbrace{-\infty,\ldots,-\infty}_{j\ \text{times}}}\{R\} = P_{\theta_1,\ldots,\theta_{k-j}}\{I^j(R)\}. \tag{44}$$

Proof of Theorem 3.1. (i) Note, for any monotone rule $E$, the smallest probability of $E_{k,j}$ over $A_j(\varepsilon)$ occurs when $\theta = w_{k,1}$ defined in (40), as well as when $\theta$ is any permutation of $w_{k,1}$. Furthermore, for any monotone rule $E$ that controls the FWER, we must have

$$P_\theta\{E_{k,j}\} \le \alpha$$

when $\theta$ is

$$v_{k,j} = (\underbrace{\infty,\ldots,\infty}_{j-1\ \text{times}},\ \underbrace{0,\ldots,0}_{k-j+1\ \text{times}}), \tag{45}$$

or permutations of $v_{k,j}$.

To prove the optimality result (30), consider another rule $E$, with $E^c_{k,1}$ the subset of $\mathbf{R}^k$ that accepts all null hypotheses. Suppose there exists $x = (x_1, \dots, x_k) \notin D^c_{k,1}$, but $x \in E^c_{k,1}$. Then there exists at least one component of $x$, say $x_1$, with $x_1 > c_{k,1}$. By monotonicity, the set

$$L(x) = \{y \in \mathbf{R}^k : y_i \le x_i\} \tag{46}$$

is also in $E^c_{k,1}$. Then

$$P_{w_{k,1}}\{E^c_{k,1}\} \ge P_{w_{k,1}}\{L(x)\} = P_\varepsilon\{X_1 \le x_1\} > P_\varepsilon\{X_1 \le c_{k,1}\} = 1 - \beta_{k,1}(\alpha, \varepsilon),$$

and so the smallest power of $E$ over $A_1(\varepsilon)$ satisfies

$$P_{w_{k,1}}\{E_{k,1}\} < \beta_{k,1}(\alpha, \varepsilon).$$

Therefore, in order for $E$ to be optimal we must have

$$E^c_{k,1} \subset D^c_{k,1}.$$

But if $D_{k,1}$ is a proper subset of $E_{k,1}$ (except for a set with 0 Lebesgue measure), then

$$P_{\underbrace{0,\dots,0}_{k \text{ times}}}\{E_{k,1}\} > P_{\underbrace{0,\dots,0}_{k \text{ times}}}\{D_{k,1}\} = \alpha,$$

a contradiction if $E$ controls the FWER. Therefore, (30) is proved.

To prove the result (31) with $j = k$, let $E$ be any other monotone decision rule which has strong control and satisfies the constraint

$$E_{k,k} \subset D_{k,k-1} = \{X_{r_1} \ge c_{k,1}, \dots, X_{r_{k-1}} \ge c_{k,k-1}\}.$$

Suppose $E_{k,k}$ includes a point $y = (y_1, \dots, y_k)$, where $y_i \ge c_{k,i}$ for $i = 1, \dots, k-1$ and $y_k < c_{k,k}$. Then

$$\begin{aligned} P_{\underbrace{\infty,\dots,\infty}_{k-1 \text{ times}},0}\{E_{k,k}\} &\ge P_{\underbrace{\infty,\dots,\infty}_{k-1 \text{ times}},0}\{X_{r_1} \ge y_1, \dots, X_{r_{k-1}} \ge y_{k-1}, X_k \ge y_k\} \\ &= P_0\{X_k \ge y_k\} > P_0\{X_k > c_{k,k}\} = \alpha, \end{aligned}$$

a contradiction of strong control. So such a point $y$ cannot be in $E_{k,k}$, nor can any permutation of the coordinates of $y$ (by a similar argument). Therefore, $E_{k,k}$ can at most include $D_{k,k}$. But taking away any points from $D_{k,k}$ could only lower the minimum power at $(\varepsilon, \dots, \varepsilon)$, and so $D_{k,k}$ is optimal.

To prove the result (31) with $1 < j < k$, let $E$ be any other monotone decision rule which has strong control and satisfies the constraint (32). Let

$$X_{j:1} \ge X_{j:2} \ge \dots \ge X_{j:j} \tag{47}$$

denote the ordered values of $X_1, \dots, X_j$. Since $E$ has strong control, it follows by (43) that

$$P_{\underbrace{0,\dots,0}_{k-j+1\ \text{times}}}\{U^{j-1}(E_{k,j})\} = \alpha.$$

Hence, $U^{j-1}(E_{k,j})$ can be viewed as a rejection region in $\mathbf{R}^{k-j+1}$ for the case with $k$ and $j$ replaced by $k' = k - j + 1$ and $j' = 1$. [Note that if $E_{k,j}$ satisfies the constraint $E_{k,j} \subset D_{k,j-1}$, then

$$U^{j-1}(E_{k,j}) \subset U^{j-1}(D_{k,j-1}) = \mathbf{R}^{k-j+1},$$

so the constraint is vacuous.] It follows that

$$P_{0,\underbrace{-\infty,\dots,-\infty}_{k-j\ \text{times}}}\{U^{j-1}(E_{k,j})\} \le \beta_{k-j+1,1}(\alpha, 0)$$

or

$$P_{0,\underbrace{-\infty,\dots,-\infty}_{k-j\ \text{times}},\underbrace{\infty,\dots,\infty}_{j-1\ \text{times}}}\{E_{k,j}\} \le \beta_{k-j+1,1}(\alpha, 0).$$

By the same reasoning applied to any permutation of

$$\theta = (0, \underbrace{-\infty, \dots, -\infty}_{k-j\ \text{times}}, \underbrace{\infty, \dots, \infty}_{j-1\ \text{times}}),$$

$$\begin{aligned} &P_{0,\underbrace{\infty,\dots,\infty}_{j-1\ \text{times}},\underbrace{-\infty,\dots,-\infty}_{k-j\ \text{times}}}\{E_{k,j}\} \\ &\quad = P_{0,\underbrace{\infty,\dots,\infty}_{j-1\ \text{times}}}\{I^{k-j}(E_{k,j})\} \le \beta_{k-j+1,1}(\alpha, 0). \end{aligned}$$

So $I^{k-j}(E_{k,j})$ is a rejection region in $\mathbf{R}^j$ that controls the Type 1 error at the point

$$(0, \underbrace{\infty, \dots, \infty}_{j-1\ \text{times}})$$

(as well as at permutations of its coordinates), not at level $\alpha$, but at level $\beta_{k-j+1,1}(\alpha, 0)$. [In words, if you use the rule $E$ which is originally designed to test $k$ hypotheses, but you ignore the last $k - j$ hypotheses, the overall probability of a Type 1 error for testing the $j$ hypotheses is reduced to $\beta_{k-j+1,1}(\alpha, 0)$.] Also, note that the constraint $E_{k,j} \subset D_{k,j-1}$ implies

$$I^{k-j}(E_{k,j}) \subset I^{k-j}(D_{k,j-1}) = \{X_{j:1} \ge c_{k,1}, \dots, X_{j:j-1} \ge c_{k,j-1}\}.$$

(Note that $c_{k,j}$ always refers to the critical values based on the given value of $\alpha$, so its dependence on $\alpha$ is suppressed.) Then, by the case with $k$ and $j$ replaced by $j$ and $j$ (already proved above) and $\alpha$ replaced by $\beta_{k-j+1,1}(\alpha,0)$, it follows that

$$P_{\underbrace{\varepsilon,\dots,\varepsilon}_{j\text{ times}}}\{I^{k-j}(E_{k,j})\} \le \beta_{j,j}(\beta_{k-j+1,1}(\alpha,0),\varepsilon)$$

or

$$P_{\underbrace{\varepsilon,\dots,\varepsilon}_{j\text{ times}},\underbrace{-\infty,\dots,-\infty}_{k-j\text{ times}}}\{E_{k,j}\} \le \beta_{j,j}(\beta_{k-j+1,1}(\alpha,0),\varepsilon). \tag{48}$$

We must argue that the right-hand side of (48) is $\beta_{k,j}(\alpha,\varepsilon)$. But notice that if we apply the above reasoning to $E = D$, the inequalities are all equalities. Indeed,

$$U^{j-1}(D_{k,j}) = \{\text{at least one of } X_1,\dots,X_{k-j+1} \ge c_{k,j}\}$$

and the optimal minimum power (with $\varepsilon = 0$) for the subproblem with $k' = k - j + 1$, $j' = 1$ and $\alpha' = \alpha$ is $\beta_{k-j+1,1}(\alpha,0)$. Also,

$$I^{k-j}(D_{k,j}) = \{X_{j:1} \ge c_{k,1},\dots,X_{j:j} \ge c_{k,j}\}$$

is optimal for the case $k'' = j'' = j$ at the level $\alpha'' = \beta_{k-j+1,1}(\alpha,0)$. Indeed, checking the level condition,

$$\begin{aligned} P_{0,\underbrace{\infty,\dots,\infty}_{j-1\text{ times}}}\{I^{k-j}(D_{k,j})\} &= P_0\{X_1 \ge c_{k,j}\} \\ &= P_0\{X_1 \ge c_{k-j+1,1}\} = \beta_{k-j+1,1}(\alpha,0). \end{aligned}$$

So, by the case $k'' = j'' = j$,

$$\begin{aligned} P_{\underbrace{\varepsilon,\dots,\varepsilon}_{j\text{ times}}}\{I^{k-j}(E_{k,j})\} &\le P_{\underbrace{\varepsilon,\dots,\varepsilon}_{j\text{ times}}}\{I^{k-j}(D_{k,j})\} \\ &= P_{\underbrace{\varepsilon,\dots,\varepsilon}_{j\text{ times}}}\{X_{j:1} \ge c_{k,1},\dots,X_{j:j} \ge c_{k,j}\} = \beta_{k,j}(\alpha,\varepsilon). \end{aligned}$$

The proof of (ii) is completely analogous to the proof of Theorem 2.1(iii), with the help of Lemma 3.1(ii). □

PROOF OF LEMMA 4.1. To prove (i), suppose $H_1,\dots,H_p$ are true and $H_{p+1},\dots,H_k$ are false. A Type 1 error occurs if any of $H_1,\dots,H_p$ are rejected. For the rule $D$, the set where any of $H_1,\dots,H_p$ is rejected is a monotone increasing set (invoking the monotonicity of critical values). Hence the probability of this event is largest under this configuration of true and false hypotheses when

$$(\theta_1,\dots,\theta_k) = (\underbrace{0,\dots,0}_{p\text{ times}},\underbrace{\infty,\dots,\infty}_{k-p\text{ times}}),$$

and this probability is equal to

$$
\begin{aligned}
&P_{\underbrace{0,\ldots,0}_{p\ \text{times}},\underbrace{\infty,\ldots,\infty}_{k-p\ \text{times}}}\{\text{reject any of } H_1,\ldots,H_p\}\\
&\quad= P_{\underbrace{0,\ldots,0}_{p\ \text{times}},\underbrace{\infty,\ldots,\infty}_{k-p\ \text{times}}}\{\text{reject any of } H_1,\ldots,H_p \\
&\qquad\qquad \cap \text{reject all of } H_{p+1},\ldots,H_k\},
\end{aligned}
\tag{49}
$$

because as $(\theta_1,\ldots,\theta_k)\to(\underbrace{0,\ldots,0}_{p\ \text{times}})$, $X_{(p+1)} > d_{p+1}$ with probability tending to one, and so the hypotheses $H_{p+1},\ldots,H_k$ are rejected with probability tending to one. Then (49) is bounded above by

$$
\begin{aligned}
&P_{\underbrace{0,\ldots,0}_{p\ \text{times}},\underbrace{\infty,\ldots,\infty}_{k-p\ \text{times}}}\{\text{at least } k-p+1 \text{ rejections}\}\\
&\quad= P_{\underbrace{0,\ldots,0}_{p\ \text{times}},\underbrace{\infty,\ldots,\infty}_{k-p\ \text{times}}}\{D_{k,k-p+1}\}\\
&\quad= 1-P_{\underbrace{0,\ldots,0}_{p\ \text{times}},\underbrace{\infty,\ldots,\infty}_{k-p\ \text{times}}}\{X_{(1)}\le d_1,\ldots,X_{(p)}\le d_p\}\\
&\quad= 1-P_{\underbrace{0,\ldots,0}_{p\ \text{times}}}\{L_j\}=\alpha,
\end{aligned}
$$

by (34) and (35).

To prove (ii), note that the minimum power occurs when $\theta$ is one of the $\binom{k}{j}$ points with $j$ values of $\varepsilon$ and $k-j$ values of $-\infty$, such as $w_{k,j}$ given by (40). Then

$$
\begin{aligned}
P_{w_{k,j}}(D_{k,j}) &= 1-P_{w_{k,j}}\{X_{(1)}\le d_1,\ldots,X_{(k-j+1)}\le d_{k-j+1}\}\\
&= 1-P_{\underbrace{\varepsilon,\ldots,\varepsilon}_{j\ \text{times}}}\{\{X_1\le d_{k-j+1}\}\cup\cdots\cup\{X_j\le d_{k-j+1}\}\},
\end{aligned}
$$

which reduces to $\tilde{\beta}_{k,j}(\alpha,\varepsilon)$ as claimed. $\square$

Proof of Theorem 4.1. To prove (37) (the case $j=k$), first observe that

$$D_{k,k}=\{X_{(1)}>d_1\}.$$

Consider another monotone rule $E$, and suppose there exists some point $x=(x_1,\ldots,x_k)$ with $x\in E_{k,k}$ but $x\notin D_{k,k}$. Then there exists at least one component of $x$, say $x_1$, with $x_1<d_1$. By monotonicity the set

$$M(x)=\{y\in \mathbf{R}^k : y_i\ge x_i\}$$

is also in $E_{k,k}$. Then

$$P_{0,\underbrace{\infty,\dots,\infty}_{k-1\ \text{times}}}\{E_{k,k}\} \geq P_{0,\underbrace{\infty,\dots,\infty}_{k-1\ \text{times}}}\{M(x)\}$$

$$= P_0\{X_1 \geq x_1\} > P_0\{X_1 \geq d_1\} = \alpha,$$

which would contradict strong control. So we must have $E_{k,k} \subset D_{k,k}$. But then

$$P_{\underbrace{\varepsilon,\dots,\varepsilon}_{k\ \text{times}}}\{E_{k,k}\} \leq P_{\underbrace{\varepsilon,\dots,\varepsilon}_{k\ \text{times}}}\{D_{k,k}\},$$

and so (37) is proved.

To prove the result (38) in the case $j = 1$, the constraint is that $E_{k,1}$ must contain $D_{k,2}$, or, equivalently,

$$E_{k,1}^c \subset D_{k,2}^c = \bigcap_{i=1}^{k-1} \{X_{(i)} \leq d_i\}.$$

Suppose $x = (x_1, \dots, x_k) \in E_{k,1}^c$ but $x \notin D_{k,1}^c$. For the sake of argument, assume the $x_i$ are nondecreasing in $i$ with $x_i \leq d_i$ for $i = 1, \dots, k-1$ (so the constraint is satisfied), but $x_k > d_k$. Then $x \in E_{k,1}^c$ implies $L(x) \in E_{k,1}^c$, where $L(x)$ is defined in (46). So

$$P_{\underbrace{-\infty,\dots,-\infty}_{k-1\ \text{times}},\varepsilon}\{E_{k,1}^c\} \geq P_{\underbrace{-\infty,\dots,-\infty}_{k-1\ \text{times}},\varepsilon}\{L(x)\}$$

$$= P_\varepsilon\{X_k \leq x_k\} > P_\varepsilon\{X_k \leq d_k\}.$$

Therefore

$$P_{\underbrace{-\infty,\dots,-\infty}_{k-1\ \text{times}},\varepsilon}\{E_{k,1}\} < P_\varepsilon\{X_k > d_k\} = \tilde{\beta}_{k,1}(\alpha, \varepsilon),$$

and so $E_{k,1}$ is less powerful than $D_{k,1}$. Therefore such a point $x$ cannot exist in order for $E_{k,1}$ to be optimal. (A similar argument applies to any permutation of the coordinates of $x$.) Then $x \in D_{k,1}$ implies $x \in E_{k,1}$. But adding any points $x$ to $D_{k,1}$ would increase the probability of rejection when $\theta = (0, \dots, 0)$, and this would contradict the level constraint. So the case $j = 1$ is proved.

To prove (38) for $1 < j < k$, let $E$ be any other monotone decision rule which has strong control and satisfies the constraint (39). Since the set $E_{k,j}$ cannot have probability greater than $\alpha$ when $\theta = v_{k,j}$, where $v_{k,j}$ is given by (45), we must have

$$\alpha \geq P_{v_{k,j}}\{E_{k,j}\} = P_{\underbrace{0,\dots,0}_{k-j+1\ \text{times}}}\{U^{j-1}(E_{k,j})\}$$

by (43). Therefore $U^{j-1}(E_{k,j})$ is a region in $\mathbf{R}^{k-j+1}$ which has rejection probability $\alpha$ when $\theta = (0, \dots, 0)$. Note that the constraint $E_{k,j} \supset D_{k,j+1}$ implies the region $U^{j-1}(E_{k,j})$ must contain

$$U^{j-1}(D_{k,j+1}) = \{X_{k-j+1:k-j+1} > d_1\} \cup \cdots \cup \{X_{k-j+1:2} > d_{k-j}\}.$$

Therefore, by the case considered above with $k' = k - j + 1$ and $j' = 1$, the optimal region in $\mathbf{R}^{k-j+1}$ is

$$\{X_{k-j+1:k-j+1} > d_1\} \cup \cdots \cup \{X_{k-j+1:1} > d_{k-j+1}\},$$

which, in fact, is equal to $U^{j-1}(D_{k,j})$. So

$$\begin{aligned} P_{0,\underbrace{-\infty,\dots,-\infty}_{k-j \text{ times}}}\{U^{j-1}(E_{k,j})\} &\le P_{0,\underbrace{-\infty,\dots,-\infty}_{k-j \text{ times}}}\{U^{j-1}(D_{k,j})\} \\ &= P_0\{X_1 > d_{k-j+1}\} = \tilde{\beta}_{k-j+1,1}(\alpha, 0). \end{aligned}$$

Using (43) and applying the argument to any permutation of $v_{k,j}$, we have

$$P_{0,\underbrace{\infty,\dots,\infty}_{j-1 \text{ times}},\underbrace{-\infty,\dots,-\infty}_{k-j \text{ times}}}\{E_{k,j}\} \le \tilde{\beta}_{k-j+1,1}(\alpha, 0),$$

or by (44),

$$P_{0,\underbrace{\infty,\dots,\infty}_{j-1}}\{I^{k-j}(E_{k,j})\} \le \tilde{\beta}_{k-j+1,1}(\alpha, 0).$$

So $I^{k-j}(E_{k,j})$ is a rejection region in $\mathbf{R}^j$ that controls the Type 1 error at

$$(0, \underbrace{\infty, \dots, \infty}_{j-1 \text{ times}})$$

(as well as permutations of its coordinates), not at level $\alpha$, but at level $\tilde{\beta}_{k-j+1,1}(\alpha, 0)$. [Also note that the constraint $E_{k,j} \supset D_{k,j+1}$ implies

$$I^{k-j}(E_{k,j}) \supset I^{k-j}(D_{k,j+1}) = \varnothing,$$

which is always satisfied.] By the case with $k'' = j'' = j$ and $\alpha'' = \tilde{\beta}_{k-j+1,1}(\alpha, 0)$ considered above,

$$I^{k-j}(D_{k,j}) = \{\min(X_1, \dots, X_j) > d_{k-j+1}\}$$

is optimal for this case and so

$$P_{\underbrace{\varepsilon,\dots,\varepsilon}_{j \text{ times}}}\{I^{k-j}(E_{k,j})\} \le P_{\underbrace{\varepsilon,\dots,\varepsilon}_{j \text{ times}}}\{I^{k-j}(D_{k,j})\} = \tilde{\beta}_{k,j}(\alpha, \varepsilon)$$

by Lemma 4.1(ii). Therefore

$$P_{\underbrace{\varepsilon,\dots,\varepsilon}_{j \text{ times}},\underbrace{-\infty,\dots,-\infty}_{k-j \text{ times}}}\{E_{k,j}\} \le \tilde{\beta}_{k,j}(\alpha, \varepsilon),$$

as was to be proved.

The proof of (ii) is analogous to the proof of Theorem 3.1(ii). □

**Acknowledgment.** Special thanks to an anonymous referee for valuable comments and references.

E. L. Lehmann
J. P. Shaffer
Department of Statistics
University of California, Berkeley
Berkeley, California 94720
USA
e-mail: shaffer@stat.berkeley.edu

J. P. Romano
Department of Statistics
Stanford University
Stanford, California 94305-4065
USA
e-mail: romano@stat.stanford.edu